\documentclass[12pt,a4paper]{article}
\usepackage{etoolbox} \usepackage{amsfonts}
\usepackage[margin=1in]{geometry}
\usepackage{times}
\usepackage{amsthm,amssymb}
\usepackage{amsmath,xypic}
\usepackage{multicol}
\usepackage{hyphenat}
\usepackage[usenames,dvipsnames]{color}
\usepackage{rotating}
\usepackage{longtable}
\usepackage{caption}
\usepackage[most]{tcolorbox}

\usepackage{epigraph}
\usepackage{enumitem}
\setlist[enumerate]{listparindent=0.5in}
\usepackage{mathrsfs}
\DeclareMathAlphabet{\mathscrbf}{OMS}{mdugm}{b}{n} \newcommand{\be}{\begin{equation}}
\newcommand{\ee}{\end{equation}}
\newcommand{\bes}{\begin{equation*}}
\newcommand{\ees}{\end{equation*}}
\newcommand{\bea}{\begin{eqnarray}}
\newcommand{\eea}{\end{eqnarray}}
\newcommand{\beas}{\begin{eqnarray}}
\newcommand{\eeas}{\end{eqnarray}}
\newcommand{\ben}{\begin{note}}
\newcommand{\een}{\end{note}}
\newcommand{\bexl}{\vskip0.1em\noindent\hrulefill\vskip1em\begin{ExerciseList}}
\newcommand{\eexl}{\end{ExerciseList}\hrulefill}

\newcommand{\bthm}{\begin{theorem}}
\newcommand{\ethm}{\end{theorem}}
\newcommand{\bpro}{\begin{prop}}
\newcommand{\epro}{\end{prop}}
\newcommand{\bcor}{\begin{corollary}}
\newcommand{\ecor}{\end{corollary}}
\newcommand{\bcon}{\begin{conjecture}}
\newcommand{\econ}{\end{conjecture}}
\newcommand{\bp}{\begin{proof}}
\newcommand{\ep}{\end{proof}}
\newcommand{\blem}{\begin{lemma}}
\newcommand{\elem}{\end{lemma}}
\newcommand{\bn}{\begin{note}}
\newcommand{\en}{\end{note}}
\newcommand{\benum}{\begin{enumerate}}
\newcommand{\eenum}{\end{enumerate}}
\newcommand{\bed}{\begin{defn}}
\newcommand{\eed}{\end{defn}}
\newcommand{\brem}{\begin{remark}}
\newcommand{\erem}{\end{remark}}

\newcommand{\btik}{\begin{tikzpicture}\begin{axis}[scale=0.5,axis y line=center, axis x line=middle]}
\newcommand{\etik}{\end{axis}\end{tikzpicture}}

\let\mapsto=\longmapsto

\newcommand{\upperRomannumeral}[1]{\uppercase\expandafter{\romannumeral#1}}

 \usepackage[pagewise]{lineno}
\usepackage{stmaryrd}
\usepackage[usenames,dvipsnames]{color}
\usepackage[authoryear,sort&compress]{natbib}
\let\cite=\citep
\usepackage[colorlinks,citecolor=blue,colorlinks=true,hyperindex, citecolor=cyan, urlcolor=cyan]{hyperref}
\newtoggle{arxiv}
\togglefalse{arxiv}
\iftoggle{arxiv}{
}
{
}
\newtheorem{theorem}[equation]{Theorem}      \newtheorem{lemma}[equation]{Lemma}          \newtheorem{corollary}[equation]{Corollary}  \newtheorem{proposition}[equation]{Proposition}

\theoremstyle{definition}
\newtheorem{conj}[equation]{Conjecture}

\newtheorem{question}[equation]{Question}

\theoremstyle{definition}
\newtheorem{defn}[equation]{Definition}
\theoremstyle{remark}

\theoremstyle{definition}
\newtheorem{remark}[equation]{Remark}

\numberwithin{equation}{section}

\let\isom=\simeq

\let\tensor=\otimes
\newcommand{\A}{\mathcal{A}}

\newcommand{\bF}{{\bar{F}}}

\newcommand{\C}{{\mathbb C}}

\newcommand{\F}{{\mathbb F}}

\newcommand{\N}{\mathcal{N}}

\newcommand{\Q}{{\mathbb Q}}

\newcommand{\Z}{{\mathbb Z}}

\renewcommand{\int}{\operatorname{int}}
\renewcommand{\O}{{\mathcal O}}
\renewcommand{\P}{{\mathbb P}}

\newcommand{\mapright}[1]{{\xymatrix{{}\ar[r]^{#1}&{}}}}

\renewcommand{\bpro}{\begin{proposition}}
	\renewcommand{\epro}{\end{proposition}}
\renewcommand{\bcon}{\begin{conj}}
	\renewcommand{\econ}{\end{conj}}
\let\mathcal=\mathscr

\title{The Absolute Grothendieck Conjecture is false for Fargues-Fontaine Curves}
\author{Kirti Joshi}

\begin{document}
\maketitle

\begin{abstract}
	I prove that the absolute Grothendieck Conjecture is false for Fargues-Fontaine curves. 
\end{abstract}

\epigraphwidth0.55\textwidth
\epigraph{And in its time the spell was snapt: once more\\
	I viewed the ocean green,\\
	I look'd far-forth, but little saw\\ 
	Of what else had been seen.}{Samuel Taylor Coleridge \cite{coleridge}}

\newcommand{\act}{\curvearrowright}
\newcommand{\lmp}{{\Pi\act\Ot}}
\newcommand{\lmpi}{{\lmp}_{\int}}
\newcommand{\lmpf}{\lmp_F}
\newcommand{\Om}{\O^{\times\mu}}
\newcommand{\Omf}{\O^{\times\mu}_{\bF}}
\renewcommand{\N}{\mathbb{N}}
\newcommand{\yoga}{Yoga}
\newcommand{\gl}[1]{{\rm GL}(#1)}
\newcommand{\bK}{\overline{K}}
\newcommand{\reptrip}{\rho:G_K\to\gl V}
\newcommand{\reptripp}[1]{\rho\circ\alpha:G_{\ifstrempty{#1}{K}{{#1}}}\to\gl V}
\newcommand{\benumlab}{\begin{enumerate}[label={{\bf(\arabic{*})}}]}
\newcommand{\ord}{\mathop{\rm ord}\nolimits}	
\newcommand{\kcs}{K^\circledast}
\newcommand{\lcs}{L^\circledast}
\renewcommand{\A}{\mathbb{A}}
\newcommand{\bfq}{\bar{\mathbb{F}}_q}
\newcommand{\tripod}{\P^1-\{0,1728,\infty\}}

\newcommand{\vseq}[2]{{#1}_1,\ldots,{#1}_{#2}}
\newcommand{\anab}[4]{\left({#1},\{#3 \}\right)\anabelmap\left({#2},\{#4 \}\right)}

\newcommand{\gln}{{\rm GL}_n}
\newcommand{\glo}[1]{{\rm GL}_1(#1)}
\newcommand{\glt}[1]{{\rm GL_2}(#1)}

\newcommand{\iut}{\cite{mochizuki-iut1, mochizuki-iut2, mochizuki-iut3,mochizuki-iut4}}
\newcommand{\topics}{\cite{mochizuki-topics1,mochizuki-topics2,mochizuki-topics3}}

\newcommand{\linv}{\mathfrak{L}}
\newcommand{\bedef}{\begin{defn}}
\newcommand{\eedef}{\end{defn}}
\renewcommand{\act}[1][]{\overset{#1}{\curvearrowright}}
\newcommand{\bfx}{\overline{F(X)}}
\newcommand{\anabelmap}{\leftrightsquigarrow}
\newcommand{\ban}[1][G]{\mathcal{B}({#1})}
\newcommand{\pit}{\Pi^{temp}}
 
 \newcommand{\bL}{\overline{L}}
 \newcommand{\bkm}{\bK_M}
 \newcommand{\vbk}{v_{\bK}}
 \newcommand{\vbkm}{v_{\bkm}}
\newcommand{\ocs}{\O^\circledast}
\newcommand{\ot}{\O^\triangleright}
\newcommand{\ocsk}{\ocs_K}
\newcommand{\otk}{\ot_K}
\newcommand{\ok}{\O_K}
\newcommand{\oko}{\O_K^1}
\newcommand{\oks}{\ok^*}
\newcommand{\Qpb}{\overline{\Q}_p}
\newcommand{\Qpbh}{\widehat{\overline{\Q}}_p}
\newcommand{\tr}{\triangleright}
\newcommand{\ocpt}{\O_{\C_p}^\tr}
\newcommand{\ocpf}{\O_{\C_p}^\flat}
\newcommand{\sG}{\mathscr{G}}
\newcommand{\sxfe}{\mathscr{X}_{F,E}}
\newcommand{\sxfep}{\mathscr{X}_{F,E'}}
\newcommand{\loglt}{\log_{\sG}}
\newcommand{\fc}{\mathfrak{t}}
\newcommand{\ku}{K_u}
\newcommand{\kup}{\ku'}
\newcommand{\kt}{\tilde{K}}
\newcommand{\sGpf}{\sG(\O_K)^{pf}}
\newcommand{\hgm}{\widehat{\mathbb{G}}_m}
\newcommand{\bE}{\overline{E}}
\newcommand{\loccit}{{loc. cit.}}

\newcommand{\FF}{\cite{fargues-fontaine}}

\section{Introduction}
Let $E,E'$ be fields. Following \cite[Section 2]{joshi-anabelomorphy}, I say that $E$ and $E'$ are \emph{anabelomorphic} (denoted as $E\anabelmap E'$) if there exists a topological  isomorphism of their absolute Galois groups $G_E\isom G_{E'}$, and refer to a topological isomorphism $\alpha:G_E\isom G_{E'}$ as an \emph{anabelomorphism} $\alpha:E\anabelmap E'$ between $E$ and $E'$.  I will say that an anabelomorphism $E\anabelmap E'$ is a \emph{strict anabelomorphism} if $E$ is not isomorphic to $E'$. Anabelomorphism of fields is an equivalence relation and in \loccit\ the invariants of the anabelomorphism class of a field are called \emph{amphoric}, for example if $E$ is a $p$-adic field (here and elsewhere in this paper a $p$-adic field will mean a finite extension of $\Q_p$) then the residue characteristic $p$ of $E$ and the degree $[E:\Q_p]$ are amphoric (for a longer list of amphoric quantities see \cite[Theorem 2.4.3]{joshi-anabelomorphy}).

The notion of anabelomorphisms of fields can be extended to a large class of smooth schemes by replacing the absolute Galois group by the \'etale fundamental group in the definition. More precisely consider the class of geometrically connected, smooth varieties over fields.  Let $E$ be a field and $X/E$ be a geometrically connected, smooth scheme  and let $\pi_1(X/E)$ be its \'etale fundamental group (computed for some choice of geometric base-point of $X$).  I say that $X/E$ and $X'/E'$ are \emph{anabelomorphic schemes} if their \'etale fundamental groups are topologically isomorphic. 

As mentioned earlier, anabelomorphy of fields (resp. schemes) is an equivalence relation on the respective classes. Note that isomorphism of schemes is another (tautological) equivalence relation on schemes and isomorphic schemes are evidently anabelomorphic.

The extraordinary \emph{absolute Grothendieck Conjecture} (see \cite{grothendieck-anabelian}) asserts that in some situations these two equivalence relations coincide:

\begin{tcolorbox}[boxsep=10pt,bicolor,frame hidden,colback=white,boxrule=0pt]
\emph{ Any two geometrically connected, smooth hyperbolic curves over number fields are anabelomorphic if and only if they are isomorphic.} 
\end{tcolorbox}

This is known to be true and has been extended to include finite and $p$-adic fields thanks to the remarkable works of Hiroaki Nakamura \cite{nakamura90b} (the genus zero hyperbolic case),  \cite{pop94} (the  birational case), Akio Tamagawa \cite{tamagawa97-gconj}, Shinichi Mochizuki \cite{mochizuki96-gconj} (finite fields and number fields) and \cite{mochizuki99} ($p$-adic fields). The formulation  of Grothendieck's conjecture  considered above is the simplest (and is adequate for the present paper), but let me say that there are other variants of the (absolute) Grothendieck Conjecture  which are also considered in the literature on this subject and the aforementioned papers will serve as a starting point for the interested reader. 

Let me point out  that Mochizuki has established (see \cite{mochizuki99}) that the relative Grothendieck conjecture holds for smooth, hyperbolic curves over a  $p$-adic field,  and has also proved that the absolute Grothendieck Conjecture is true for smooth, hyperbolic curves of strict Belyi Type (see \cite[Corollary 2.12]{mochizuki07-cuspidalizations-proper-hyperbolic}).  Note that a hyperbolic curve of strict Belyi Type (\cite[Definition~2.9]{mochizuki07-cuspidalizations-proper-hyperbolic}) is necessarily an affine hyperbolic curve by definition, while the examples of this paper (Theorem~\ref{th:gconj}, Theorem~\ref{th:gconj-nontrivial-pi}) are complete curves (in the sense that the degree of the divisor of any non-constant meromorphic function is zero--see \cite[D\'efinition~5.1.3]{fargues-fontaine} and \cite[Th\'eor\`eme~5.2.7]{fargues-fontaine}).

In \cite[Remark 1.3.5.1]{mochizuki04} it was suggested that the absolute Grothendieck conjecture may be false for hyperbolic curves over  $p$-adic fields. However, in correspondence (2021), Mochizuki has reminded  me that \cite[Remark 1.3.5.1]{mochizuki04} was made in relation to what was known at the  time that paper was written.  Subsequent works (\cite[Corollary 2.12]{mochizuki07-cuspidalizations-proper-hyperbolic} and \cite{murotani})  have raised the possibility that the absolute Grothendieck conjecture  for hyperbolic curves over $p$-adic fields may very well be true.

\newcommand{\cdi}{\mathscr{D}^{irrat}}
\newcommand{\cdih}{\mathscr{D}^{irrat}_{hyp}}
\newcommand{\proj}{{\rm Proj}}

 Now consider the  class $\cdi$ of separated schemes $X$ satisfying the following:

\benum[label={\bf (D.\arabic{*}}),leftmargin=3\parindent]
 \item\label{hyp:1} $X$ is a Dedekind scheme i.e. one dimensional, Noetherian, and regular.
\item\label{hyp:2} $X=\proj(R)$ for a graded ring $R=\oplus_{n\geq 0}R_n$ generated by degree one elements over $R_0$, and let $\O_X(1)$ be the tautological line bundle given by this grading.
 \item\label{hyp:3} $H^1(X,\O_X(-1))\neq 0$.
\eenum

Note that $\cdi$ contains the class of smooth, proper non-rational curves over fields (hence the superscript). Further note that $\P^1\not\in \cdi$ as $H^1(\P^1,\O_{\P^1}(-1))=0$ (here $\O_{\P^1}(1)$ is the tautological line bundle of degree one given by the construction of $\P^1$) and $\cdi\supset \cdih$ where $\cdih$ is the subclass of schemes in $\cdi$ which corresponds to smooth, complete and hyperbolic curves (i.e. of genus at least two) over fields. 
 
The purpose of this note is to record the proof of  the following:
\bthm\label{th:gconj0}\hfill
\benumlab
\item\label{th:gconj0-1} The class of Fargues-Fontaine curves is contained in $\cdi$, and
\item\label{th:gconj0-2} the absolute Grothendieck conjecture is false for Fargues-Fontaine curves and hence the conjecture is false in general for the class $\cdi$.
\eenum
\ethm

In this theorem, by Fargues-Fontaine curves, I mean the ``complete'' curves constructed in \cite[Chapitre 6]{fargues-fontaine} (for more a more precise formulation see Theorem~\ref{th:gconj-2}). The curves alluded to here \emph{are not contained in $\cdih$} as they are not of finite type over their base fields, but these curves are  complete in the sense of function theory of curves: the divisor of zeros and poles of any meromorphic function  on these curves is of degree zero (see \cite[D\'efinition~5.1.3]{fargues-fontaine} and \cite[Th\'eor\`eme~5.2.7]{fargues-fontaine}). As was  established in \cite{fargues-fontaine}, these curves play a fundamental role in  $p$-adic Hodge Theory so these curves form a natural class of examples from the point of view of the theory of $p$-adic representations. Notably, in Theorem~\ref{th:gconj-nontrivial-pi}, I show that Fargues-Fontaine curves also provide examples of  strictly anabelomorphic curves (of class $\cdi$) whose \'etale fundamental group is not isomorphic to the absolute Galois group of their respective base fields. For explicit examples of Theorem~\ref{th:gconj0}, Theorem~\ref{th:gconj} and Theorem~\ref{th:gconj-nontrivial-pi} see Remark~\ref{rem:example1} and Remark~\ref{rem:example2}.

The assertion in Theorem~\ref{th:gconj0}\ref{th:gconj0-1}  is proved, amongst many other beautiful results, in \cite{fargues-fontaine} (see below for precise references). So the main result of this paper is Theorem~\ref{th:gconj0}\ref{th:gconj0-2} and this assertion will be immediate from the more precise Theorem~\ref{th:gconj} (and also Theorem~\ref{th:gconj-nontrivial-pi}) which are proved in the next section.

Given the fundamental role which $p$-adic Hodge Theory plays in Mochizuki's work on  Grothendieck's conjecture  (see \cite{mochizuki96-gconj}, \cite{mochizuki99} and his subsequent works on  related questions)  some readers may perhaps find it surprising that the absolute Grothendieck conjecture fails for the fundamental curves of $p$-adic Hodge Theory! In some sense the point   is that these curves themselves have distinguishable $p$-adic Hodge theories. Let me also say that this paper grew out of my philosophy of combining Anabelian Geometry, Perfectoid Geometry and Modern $p$-adic Hodge Theory (see \cite{joshi-formal-groups} for additional evidence in this direction). This philosophy has recently led me to a construction of arithmetic and adelic Teichmuller spaces (see \cite{joshi-teich}). As is detailed in \cite[\S 8.12]{joshi-teich-rosetta}, the main theorem of this paper provides a natural geometric genesis and meaning  for Mochiuzki's Indeterminacy Ind1 \cite[Theorem 3.11, Page 575]{mochizuki-iut3}.

I would like to thank: Peter Scholze, Laurent Fargues, Taylor Dupuy for some correspondence; Yuichiro Hoshi, Shinichi Mochizuki for some conversations around Grothendieck's Conjecture.

\section{The main theorem}\label{se:grothendieck}
Let $F$ be an algebraically closed perfectoid field of characteristic $p>0$ \cite[Definition 3.1]{scholze12-perfectoid-ihes}. Let $E$ be a $p$-adic field i.e. $E/\Q_p$ is a finite extension.  Following \cite[Definition 2.1.1]{joshi-anabelomorphy}, I say that two $p$-adic fields are \emph{anabelomorphic} if there exists a topological isomorphism $G_E\isom G_{E'}$ of their absolute Galois groups; I write this as $E\anabelmap E'$. As is observed in \cite[Definition 2.1.1(2)]{joshi-anabelomorphy}, anabelomorphism (of $p$-adic fields) is an equivalence relation on $p$-adic fields. The notion of anabelomorphism extends to schemes: two schemes are anabelomorphic if their \'etale fundamental groups are isomorphic and two anabelomorphic schemes are \emph{strictly anabelomorphic schemes} if they are anabelomorphic but  not isomorphic.

Note that there exist $p$-adic fields which are not isomorphic but are anabelomorphic (see for instance \cite{jarden79} or \cite{joshi-anabelomorphy} for examples) and hence the absolute Grothendieck conjecture is already false for $p$-adic fields; on the other hand  Mochizuki has established (see \cite{mochizuki99}) that the Grothendieck conjecture holds for smooth, hyperbolic curves over isomorphic  $p$-adic fields.  In \cite[Remark 1.3.5.1]{mochizuki04} Mochizuki has suggested that the Grothendieck conjecture may be false for hyperbolic curves over arbitrary (i.e. non-isomorphic) $p$-adic fields.

Now suppose that $E,E'$ are $p$-adic fields. Let $\sxfe$ (resp. $\sxfep$)  be the Fargues-Fontaine curve \cite[Chap 6]{fargues-fontaine} associated to $(F,E)$ and $(F,E')$ respectively.  Note that in \loccit\ this curve is denoted by $X_{F,E,\pi}$ where $\pi$ is a uniformizer for $E$ (I will suppress $\pi$ from the notation in the present paper as it is irrelevant to what is done here).
Let me remark that $\sxfe$ is not of finite type and while it is supposed to have many properties similar to $\P^1$ (see \cite[Chap 5]{fargues-fontaine}), $\sxfe$ also shares some properties of curves of genus $\geq2$. As mentioned in the Introduction, the Fargues-Fontaine curves $\sxfe$ are  complete curves in the sense of \cite[Definition 5.1.3 and Theorem 5.2.7]{fargues-fontaine} i.e. the divisor of any meromorphic function on $\sxfe$ has degree zero.

\bp[Proof of Theorem~\ref{th:gconj0}]
Let me note that $\sxfe\in\cdi$  for every $p$-adic field and every algebraically closed perfectoid field $F$. This is proved in \cite{fargues-fontaine}: by \cite[Definition 5.1.1 and Theorem 6.5.2]{fargues-fontaine} $\sxfe$ satisfies \ref{hyp:1}; by \cite[Definition 6.1.1 and Theorem~6.5.2]{fargues-fontaine}, $\sxfe$ satisfies \ref{hyp:2}. That \ref{hyp:3} holds follows from the computation of the cohomology of the tautological line bundle $\O_{\sxfe}(1)$ on  $\sxfe$ is computed in \cite[Section 8.2.1.1]{fargues-fontaine}. This proves the assertion Theorem~\ref{th:gconj0}\ref{th:gconj0-1}. The assertion Theorem~\ref{th:gconj0}\ref{th:gconj0-2} is evident from Theorem~\ref{th:gconj} proved below.
\ep

The main theorem is the following:
\bthm\label{th:gconj} 
Assume $F$ is an algebraically closed perfectoid field of characteristic $p>0$, $E,E'$ are $p$-adic fields. Let $\alpha: E'\anabelmap E$ be an  anabelomorphism (i.e. one has an isomorphism $\alpha:G_{E'}\mapright{\isom} G_E$ of topological groups). Then one has the following assertions.
\benumlab
\item\label{th:gconj-1} There is an isomorphism of topological groups $$\pi_1(\sxfe/E)\isom G_E\isom G_{E'}\isom \pi_1(\sxfep/E').$$
\item\label{th:gconj-2} Hence $\sxfe/E$ and $\sxfep/E'$ are anabelomorphic, one dimensional Dedekind schemes over anabelomorphic $p$-adic fields $E\anabelmap E'$.
\item\label{th:gconj-3} If $E'\anabelmap E$ is a strict anabelomorphism  (i.e. $E'$ is not isomorphic to $E$) then $\sxfe$ and $\sxfep$ are not isomorphic as schemes. 
\item\label{th:gconj-4} In particular  the absolute Grothendieck Conjecture is false for Fargues-Fontaine curves in general.
\eenum
\ethm
\bp 
The first and the second assertion follows from the computation of the fundamental group of $\sxfe$ and $\sxfep$  (for $F$ algebraically closed) in \cite{fargues-fontaine}, \cite[Prop.~5.2.1]{fargues-fontaine-beijing}. So it remains to prove the third assertion (which obviously implies the fourth assertion). So let me prove the third assertion.

I provide two different proofs of this.

Suppose $$\alpha: \pi_1(\sxfe/E)\isom \pi_1(\sxfep/E')$$  is an anabelomorphism of $\sxfe/E\anabelmap \sxfep/E'$. By the identification of $\pi_1(\sxfe/E)\isom G_E$ one sees that $\alpha$ induces an isomorphism $\alpha: G_{E'} \to G_E$ hence the  fields $E'$ and $E$ are anabelomorphic.
 
Let me note a useful consequence of the fact that one has, in the present case, an anabelomorphism $E\anabelmap E'$. Let $E\supseteq E_0$ (resp. $E\supseteq E_0'$) be the maximal unramified subextensions of $E$ (resp. $E'$). Then the two  extensions $E_0,E_0'$ of $\Q_p$ are isomorphic $E_0\isom E_0'$. This is because  there is a unique unramified extension of $\Q_p$ of a given degree and as $E\anabelmap E'$ by \cite{jarden79} the degree of the maximal unramified subextensions is amphoric (i.e. determined by the topological group $G_{E'}\isom G_{E}$).  
 
Now returning to the proof of the assertion, assume that the Grothendieck conjecture is  true in this context: this means the anabelomorphism 
$$\alpha: \pi_1(\sxfe/E)\isom \pi_1(\sxfep/E')$$  induces an isomorphism of schemes $$\alpha:\sxfep \mapright{\isom} \sxfe.$$
 
By \cite{fargues-fontaine} one has $H^0(\sxfe,\O_{\sxfe})=E$. This is a part of the more general assertion  (see \cite[Chap 8, 8.2.1.1]{fargues-fontaine}) that the graded ring $P=\bigoplus_{d\geq0}P_d$ is identified with the graded ring $$P=\bigoplus_{d\in\Z}H^0(\sxfe,\O_{\sxfe}(d)), $$
with $P_d=H^0(\sxfe,\O_{\sxfe}(d))$.

Thus the isomorphism $\sxfe\isom \sxfep$ of schemes provides an isomorphism $$H^0(\sxfe,\O_{\sxfe})\isom H^0(\sxfep,\O_{\sxfep}).$$ Hence this gives us an isomorphism of rings
$$E\isom H^0(\sxfe,\O_{\sxfe})\isom H^0(\sxfep,\O_{\sxfep}) \isom E',$$
and this evidently extends to an isomorphism of these fields and by \cite{schmidt33} (\cite{lang-algebra}) any (arbitrary) isomorphism of fields equipped with a discrete valuations and complete with respect to the respective valuation topologies, is in fact an isomorphism of discretely valued fields. On the other hand I have assumed in my hypothesis \ref{th:gconj-3} that the anabelomorphism $E\anabelmap E'$ is strict i.e. $E$ is not isomorphic to $E'$ and so one has arrived at a contradiction.
\ep 

\bp[Second Proof of Theorem~\ref{th:gconj}]
Let me provide a second more natural proof which illustrates precisely how  $p$-adic Galois representations  are responsible for the failure of Grothendieck conjecture for Fargues-Fontaine curves.
 
\newcommand{\sV}{\mathcal{V}}
The idea is to use (1) on one hand the correspondence  established by Fargues-Fontaine in \cite[Chap 11]{fargues-fontaine}, \cite{fargues-fontaine-beijing}  between de Rham (resp. semi-stable and crystalline) representations $\rho:G_E \to GL(V)$ (with $V/\Q_p$ a finite dimensional vector space) of $G_E$ and  $G_E$-equivariant vector bundles of a suitable sort on $\sxfe$ (see \cite[Chap 8]{fargues-fontaine} for details).  This correspondence is given by $V\mapsto \sV=V\tensor \O_{\sxfe}$ and the Fontaine functor $D_{cris}(V)$ is naturally identified as 
$$H^0(\sxfe,\sV)= D_{cris}(V)\tensor_{E_0} E,$$
see \cite[Chap 11]{fargues-fontaine}, \cite[Theorem 6.3]{fargues-fontaine-beijing}.
By \cite{fontaine94c}, $V$ is crystalline if and only if $$\dim_{\Q_p}(V)=\dim_{E_0} D_{cris}(V).$$

(2) on the other hand a fundamental fact implicit  in the proof of \cite[Corollary 3.4, Remark~3.3.1]{hoshi13} and  \cite[Discussion on Page 3]{hoshi18} shows that if $\alpha: E'\anabelmap E$ is a strict anabelomorphism then there exist a potentially crystalline representation $\rho$ of $G_E$ such that $\rho'=\rho\circ\alpha$ is not Hodge-Tate representation of $G_{E'}$ (let me note that in Hoshi's proof,  the potentially crystalline $\Q_p$-representation, over a suitable open subgroup, is  the crystalline $\Q_p$-representation arising from a Lubin-Tate group over  $E$).
Let me give a proof now assuming that this representation is in fact crystalline (other wise one can pass to a finite extensions of $E$ over which this happens and replacing $E'$ by a suitable finite extension (denoted again by $E,E'$) such that $E'\anabelmap E$, $\rho$ is crystalline and $\rho'=\rho\circ\alpha$ is not Hodge-Tate). Choose such a crystalline representation $\rho$ of $G_E$.

Now the pull-back of $\sV$ by the isomorphism $\alpha:\sxfep\isom\sxfe$, denoted $\sV'=\alpha^*(\sV)$ (with $V'$ for the underlying vector space of the corresponding representation), evidently satisfies
$$H^0(\sxfep,\sV')\isom H^0(\sxfe,\sV).$$
Now from the identification $H^0(\sxfe,\sV)\isom D_{cris}(V)\tensor_{E_0} E$, and the fact that $E\anabelmap E'$ one knows (from \cite{jarden79} or \cite[Theorem~3.3]{joshi-anabelomorphy}) that $E_0\isom E_0'$ and also $[E:\Q_p]=[E':\Q_p]$ (i.e. anabelomorphic $p$-adic fields have the same degree over $\Q_p$) and also $[E:E_0]=[E':E_0']$ (i.e. anabelomorphic $p$-adic fields have the same absolute ramification index). So the identification of the two cohomologies gives an equality of dimensions
$$ \dim_{E_0} D_{cris}(V)\cdot [E:E_0]=\dim_{E_0'} D_{cris}(V')\cdot [E':E_0'],$$ hence one sees that $\sV'$ is also crystalline as $\dim_{\Q_p}(V)$ is the common dimension (over $E_0\isom E_0'$) of both of these vector spaces.

By the functoriality of the constructions of \cite[Chap 11]{fargues-fontaine}, the bundle $\sV'$ is the bundle corresponding to the pull-back via $\alpha$ of the representation $\rho$ of $G_E$ to $G_{E'}$ i.e to the $G_{E'}$ representation $\rho'$. Hence $\rho'$  is crystalline and hence de Rham and hence Hodge-Tate. This contradicts  my assumption that $\rho'$ is not Hodge-Tate.
\ep

The hypothesis in Theorem~\ref{th:gconj} that $F$ is an algebraically closed perfectoid field of characteristic $p>0$ can be replaced by the weaker assumption that $F$ is a perfectoid field of characteristic $p>0$. The same proof as above also proves this general case:

\brem\label{rem:example1}
To make Theorem~\ref{th:gconj} completely explicit, one can take $$F=\widehat{\overline{\F_p((t))}}$$ i.e. one can take $F$ to be the completion of the algebraic closure of $\F_p((t))$ (this field, of characteristic $p>0$, is a complete valued field and also perfect and hence perfectoid \cite{scholze12-perfectoid-ihes}). One can take for $p>2$ $$E=\Q_p(\zeta_p,\sqrt[p]{p})\text{ and } E'=\Q_p(\zeta_p,\sqrt[p]{1+p}),$$ where $\zeta_p$ is a primitive $p^{th}$-root of unity; and for $p=2$ one can take
$$E=\Q_2(\zeta_8,\sqrt{\zeta_4-1})\text{ and } E'=\Q_2(\zeta_4,\sqrt[4]{2}),$$
where $\zeta_8$ (resp. $\zeta_4$) is a primitive $8^{th}$-root of unity (resp. a primitive $4^{th}$-root of unity). Then  $E$ and $E'$ are strictly anabelomorphic $p$-adic fields. These examples are due to \cite{jarden79}--similar examples were also found in \cite{yamagata76}; for infinitely  many examples (for each prime $p$) see \cite[Lemma~4.4]{joshi-anabelomorphy}.
\erem

\bthm\label{th:gconj-nontrivial-pi}
Let $F$ be a perfectoid field of characteristic $p>0$. Let $G_F$ be the absolute Galois group of $F$.  Let $E,E'$ be $p$-adic fields. If $E\anabelmap E'$ is a strict anabelomorphism of  $p$-adic fields then  $\sxfe,\sxfep$ are anabelomorphic schemes of class $\cdi$ with  $$\pi_1(\sxfe) \isom G_F\times G_E\isom G_{F}\times G_{E'}\isom \pi_1(\sxfep)$$  but $\sxfe,\sxfep$ are not isomorphic as schemes.
\ethm

\bp 
The proof is the same as the one given above except for the assertion about the fundamental group $\pi_1(\sxfe)\isom G_F\times G_E$ which   can be found in  \cite[Prop.~5.2.1]{fargues-fontaine-beijing}. Now  by \cite[Chap 7, 7.2]{fargues-fontaine} one has $H^0(\sxfe,\O_{\sxfe})=E$. 

Thus any isomorphism $\sxfe\isom \sxfep$ provides an isomorphism of rings
$$E\isom H^0(\sxfe,\O_{\sxfe})\isom H^0(\sxfep,\O_{\sxfep}) \isom E'.$$
This extends to an isomorphism of (discretely valued) fields $E\isom E'$ as before and this contradicts my assumption that $E\anabelmap E'$ is a strict anabelomorphism.
\ep

\brem\label{rem:example2} To make   Theorem~\ref{th:gconj-nontrivial-pi} completely explicit, one can take $F$ to be the completion (with respect to its valuation) of the perfection $\F_p((t))^{\rm perf}$ of $\F_p((t))$ and $E,E'$ as in Remark~\ref{rem:example1}.
\erem
\brem
Theorem~\ref{th:gconj-nontrivial-pi} provides counter examples to the absolute Grothendieck conjecture in which the \'etale fundamental group is strictly larger than the absolute Galois groups of the base field.
\erem
\newcommand{\sxqp}{\mathscr{X}_{F,\Q_p}}
\newcommand{\Isom}{{\rm Isom}}
\brem\label{rem:general-compact-case}
Let me remark that the proof of Theorem~\ref{th:gconj} also provides a proof of the following assertion. Assume $X/E,Y/E'$ are geometrically connected, smooth, proper, hyperbolic curves over $p$-adic fields. Suppose that the absolute Grothendieck conjecture is true for $X/E,Y/E'$ i.e. one has a bijection of sets
$$\Isom_{{\Z-sch}}(X,Y)\isom \Isom^{Out}(\pi_1({X/E}),\pi_1({Y/E'})).$$
Then one has the following dichotomy:
\benumlab
\item Either one (and hence both) of the above $\Isom$-sets is an empty set i.e. the absolute Grothendieck Conjecture holds vacuously for $X/E,Y/E'$, or
\item one has an isomorphism  $E\isom E'$ of discretely valued fields.
\eenum
In particular it follows that if $E,E'$ are strictly anabelomorphic $p$-adic fields, then for geometrically connected, smooth, proper, hyperbolic curves over $E$ (resp. $E'$), the absolute Grothendieck conjecture either holds vacuously or it does not hold i.e. there exist some  (compact) hyperbolic curves  $X/E, Y/E'$  such that $$\Isom^{Out}(\pi_1({X/E}),\pi_1({Y/E'}))\neq \emptyset$$ while $$\Isom_{{\Z-sch}}(X,Y)=\emptyset.$$
\erem
\brem\label{rem:general-FF-case} Let $E,E'$ be strictly anabelomorphic $p$-adic fields (see Remark~\ref{rem:example1}). For any  geometrically connected, smooth, projective  scheme $X/\Q_p$, write $X_E=X\times_{\Q_p}E$ (resp. $X_{E'}=X\times_{\Q_p}E')$. Now let me note that the examples of Theorem~\ref{th:gconj} arise as follows. Let $\sxqp$ be the complete Fargues-Fontaine curve for the datum $F,\Q_p$.  Let $\sxfe,\sxfep$ be as in Theorem~\ref{th:gconj}. Then by  \cite[Th\'eor\`eme 6.5.2(2)]{fargues-fontaine}, one has isomorphisms of $E$-schemes (resp. $E'$-schemes) 
\be
	\label{eq:yoked-examples}\sxfe\isom \sxqp\times_{\Q_p}E \text{ and }   \sxfep\isom \sxqp\times_{\Q_p}E',
\ee
and by Theorem~\ref{th:gconj}, a strict anabelomorphism $\pi_1(\sxfe)\isom \pi_1(\sxfep)$. Hence one has
\be 
\Isom^{Out}(\pi_1({\sxfe}),\pi_1({\sxfep}))\not=\emptyset.
\ee
On the other hand by  Theorem~\ref{th:gconj}  
one has 
$$\Isom_{{\Z-sch}}(\sxfe,\sxfep)=\emptyset.$$

Remark~\ref{rem:general-compact-case} and the examples of Theorem~\ref{th:gconj} thus suggests the following  question:
\begin{question}\label{question}
Is it true that for every pair of strictly anabelomorphic $p$-adic fields $E,E'$ and for every geometrically connected, smooth, proper, hyperbolic curve $X/\Q_p$ one has 
\be 
\Isom^{Out}(\pi_1({X_E}),\pi_1({X_{E'}}))=\emptyset?
\ee
\end{question}
If the answer to Question~\ref{question} is false for some hyperbolic curve $X/\Q_p$ (or over a finite extension of $\Q_p$) then one has arrived at examples of proper hyperbolic curves of finite type providing counter examples to the absolute Grothendieck conjecture (by the method of proof of \ref{th:gconj}). 
\erem

\bibliographystyle{plainnat}

\end{document}